\documentclass{amsart}
\usepackage{amsmath}
\usepackage{amsthm}
\usepackage{amssymb,latexsym}
\newtheorem{theorem}    {Theorem}       [section]
\newtheorem{lemma}      [theorem]       {Lemma}
\newtheorem{corollary}  [theorem]       {Corollary}
\newtheorem{proposition}[theorem]       {Proposition}
\newtheorem{example}    [theorem]       {Example}

\newtheorem{definition} [theorem]       {Definition}

\def\A{{\mathbb A}}
\def\N{{\mathbb N}}
\def\Z{{\mathbb Z}}
\def\F{{\mathbb F}}
\def\Q{{\mathbb Q}}
\def\R{{\mathbb R}}

\def\Re{\operatorname{Re}}
\def\re{\Re}
\def\Gl{\operatorname{Gl}}
\def\tr{\operatorname{tr}}
\def\vol{\operatorname{vol}}
\def\dim{\operatorname{dim}}
\def\res{\operatorname{res}}

\def\Sl{\operatorname{Sl}}
\def\HNF{\operatorname{HNF}}
\begin{document}

\title{Volumes of symmetric spaces via lattice points}

\thanks{Supported by NSF grants DMS 01-00587 and 99-70085.}

\author[Gillet]{Henri Gillet}
\address{University of Illinois at Chicago}
\urladdr{http://www.math.uic.edu/\char`\~henri/}
\email{henri\char`\@math.uic.edu}

\author[Grayson]{Daniel R. Grayson}
\address{University of Illinois at Urbana-Champaign}
\urladdr{http://www.math.uiuc.edu/\char`\~dan/}
\email{dan\char`\@math.uiuc.edu}


\date{\ifcase\month\or January\or February\or March\or April\or May\or
  June\or July\or August\or September\or October\or November\or December\fi
  \space\number\day, \number\year}

\maketitle

\section*{Introduction}

In this paper we show how to use elementary methods to prove that the volume
of $\Sl_k\R/\Sl_k\Z$ is $\zeta(2) \zeta(3) \cdots \zeta(k) / k$; see Corollary \ref{volzeta}.
Using a version of reduction theory presented in this paper, we can compute 
the volumes of certain unbounded regions in Euclidean space by counting
lattice points and then appeal to the machinery of Dirichlet series to get
estimates of the growth rate of the number of lattice points appearing in the
region as the lattice spacing decreases.

In section \ref{padicvol} we present a proof of the closely related result that the Tamagawa
number of $\Sl_{k,\Q}$ is $1$ that is somewhat simpler and more
arithmetic than Weil's in \cite{MR83m:10032}.  His proof proceeds by
induction on $k$ and appeals to the Poisson summation formula, whereas the
proof here brings to the forefront local versions (\ref{localTam}) of the
formula, one for each prime $p$, which help to illuminate the appearance of
values of zeta functions in formulas for volumes.

The volume computation above is known; see, for example, \cite{Siegel36}, and
formula (24) in \cite{siegel45}.  
The methods used in the computation of the volume of $\Sl_k\R/\Sl_k\Z$ in the 
book \cite[Lecture XV]{MR91d:11070} have a different flavor from ours and do not involve
counting lattice points.  One positive point about the proof there is that it proceeds by
induction on $k$, making clear how the factor $\zeta(k)$ enters in at $k$-th
stage.
See also \cite[\S 14.12, formula
(2)]{MR99g:20090}.  The proof offered there seems to have a gap which consists of
assuming that a certain region (denoted by $T$ there) is bounded, thereby allowing
the application of \cite[\S 14.4, Theorem 3]{MR99g:20090}\footnote{called
Dirichlet's Principle in \cite[\S 5.1, Theorem 3]{MR33:4001}}.  The region in
Example \ref{unboundedexample} below shows that filling the gap is not easy, hence if
we want to compute the volume by counting lattice points, something like our use
of reduction theory in Section \ref{sec:redthy} is needed.

An almost equivalent result was proved by Minkowski --- he computed the volume
of $SO(k) {\setminus} \Sl_k\R/\Sl_k\Z$.  The relationship between the two volume
computations is made clear in the proof of \cite[\S 14.12, Theorem
2]{MR99g:20090}.

Some of the techniques we use were known to Siegel, who used similar methods in his
investigation of representability of integers by quadratic forms in
\cite{siegel37,siegel38,siegel39}.  See especially \cite[Hilfssatz 6,
p.~242]{siegel38}, which is analogous to our Lemma \ref{cusps} and the reduction
theory of Section
\ref{sec:redthy}, where we show how to compute the volume of certain unbounded
domains in Euclidean space by counting lattice points; see also the
computations in \cite[\S 9]{siegel37}, which have the same general
flavor as ours.  See also
\cite[p.~581]{MR6:38b} where Siegel omits the laborious study, using reduction
theory, of points at infinity; it is those details that concern us here.

We thank Harold Diamond for useful information about Dirichlet series and Ulf
Rehmann for useful suggestions, advice related to Tamagawa numbers, and
clarifications of Siegel's work.

\section{Counting with zeta functions}

As in \cite{MR89m:11084} we define the zeta function of a group $G$ by
summing over the subgroups $H$ in $G$ of finite index.
\begin{equation}\label{zetadef}
        \zeta(G,s) = \sum_{H \subseteq G} [G:H]^{-s}
\end{equation}
Evidently, $\zeta(\Z,s) = \zeta(s)$ and the series converges for $s > 1$.
For good groups $G$ the number of subgroups of index at most $T$ grows slowly
enough as a function of $T$ that $\zeta(G,s)$ will converge for $s$
sufficiently large.

Let's pick $k \ge 0$ and compute $\zeta(\Z^k,s)$.  Any subgroup $H$ of $\Z^k$
of finite index is isomorphic to $\Z^k$; choosing such an isomorphism amounts
to finding a matrix $A : \Z^k \to \Z^k$ whose determinant is nonzero and
whose image is $H$.  Any two matrices $A$, $A'$ with the same image $H$ are
related by an equation $A' = A S$ where $S \in \Gl_k\Z$.  

Thus the terms in the sum defining $\zeta(\Z^k,s)$ correspond to the orbits
for the action of $\Gl_k\Z$ via column operations on the set of $k\times
k$-matrices with integer entries and nonzero determinant.  A unique
representative from each orbit is provided by the matrices $A$ that are in
{\sl Hermite normal form} (see \cite[p.~66]{MR94i:11105} or
\cite[II.6]{MR49:5038}), i.e., those matrices $A$ with $A_{ij} = 0$ for $i >
j$, $A_{ii} > 0$ for all $i$, and $0 \le A_{ij} < A_{ii}$ for $ i < j $.

Let $\HNF$ be the set of integer $k\times k$ matrices in Hermite normal form.
Given positive integers $n_1, \dots, n_k$, consider the set of
matrices $A$ in $\HNF$ with $A_{ii} = n_i$ for all $i$.  The number of matrices
in it is $n_1^{k-1}  n_2^{k-2}  \cdots  n_{k-1}^1 n_k^0$.  Using that, we
compute formally as follows.
\begin{equation}\label{zetalatt}
 \begin{split}
        \zeta(\Z^k,s)
                & = \sum_{H \subseteq \Z^k} [\Z^k:H]^{-s} \\
                & = \sum_{A \in \HNF} (\det A)^{-s} \\
                & = \sum_{n_1>0, \dots, n_k>0} 
                        (n_1^{k-1}  n_2^{k-2}  \cdots  n_{k-1}^1 n_k^0)
                         (n_1 \cdots n_k)^{-s} \\
                & = \sum_{n_1>0, \dots, n_k>0} 
                        n_1^{k-1-s}  
                        n_2^{k-2-s}  \cdots 
                        n_{k-1}^{1-s} 
                        n_k^{-s} \\
                & =     \sum_{n_1>0} n_1^{k-1-s} 
                        \sum_{n_2>0} n_2^{k-2-s}  \cdots 
                        \sum_{n_{k-1}>0} n_{k-1}^{1-s} 
                        \sum_{n_{k}>0} n_{k}^{-s} \\
                & = \zeta(s-k+1)  \zeta(s-k+2)  
                        \cdots  \zeta(s-1)  \zeta(s)
 \end{split}
\end{equation}
The result $\zeta(s-k+1)  \zeta(s-k+2) \cdots  \zeta(s-1)  \zeta(s)$  is a
product of Dirichlet series with positive coefficients that converge for
$s>k$, and thus $\zeta(\Z^k,s)$ also converges for $s>k$.
This computation is old, and appears in various guises.  See, for example: 
proof 2 of Proposition 1.1 in \cite{MR89m:11084};
Lemma 10 in \cite{MR21:1966};
formula (1.1) in \cite{MR57:286};
page 64 in \cite{MR47:3318};
formula (5) and the lines following it in \cite{Siegel36}, where the
counting argument is attributed to Eisenstein, and its generalization to
number rings is attributed to Hurwitz;
and
pages 37--38 in \cite{MR83m:10032}.

\begin{lemma}\label{HNFcount}
  $
        \# \{ H \subseteq \Z^k \mid [\Z^k:H] \le T \} 
        \sim
        \zeta(2)  \zeta(3)  \cdots  \zeta(k) T^k / k
  $ for $k \ge 1$.
\end{lemma}

The right hand side is interpreted as $T$ when $k = 1$.
The notation $f(T) \sim g(T)$ means that $\lim_{T \to \infty}
f(T)/g(T) = 1$.

\begin{proof}
  We give two proofs.
  
  The first one is more elementary, and was told to us by Harold Diamond.  Writing $\zeta(s-k+1) =
  \sum n^{k-1} n^{-s}$ and letting $B(T) = \sum_{n \le T} n^{k-1}$ be the corresponding coefficient
  summatory function we see that $B(T) = T^k/k + O(T^{k-1})$.  
  If $k \ge 3$ we may apply Theorem \ref{product} to show that the coefficient summatory function
  for the Dirichlet series $\zeta(s) \zeta(s-k+1)$ behaves as $\zeta(k) T^k/k + O(T^{k-1})$.
  Applying it several more times shows that the coefficient summatory function for the Dirichlet
  series $\zeta(s) \zeta(s-1) 
  \cdots \zeta(s-k+3) \zeta(s-k+1)$ behaves as $\zeta(k) \zeta(k-1) \cdots \zeta(3) T^k/k +
  O(T^{k-1})$.  Applying it one more time we see that the coefficient summatory function for
  $\zeta(\Z^k,s) = \zeta(s) \cdots \zeta(s-k+1)$ behaves as $\zeta(k) \zeta(k-1) \cdots \zeta(2)
  T^k/k + O(T^{k-1} \log T)$, which in turn implies the result.
  
  The second proof is less elementary, since it uses a Tauberian theorem.  From (\ref{zetalatt}) we
  know that the rightmost (simple) pole of $\zeta(\Z^k,s)$ occurs at $s = k$, that the residue
  there is the product $ \zeta(2) \zeta(3) \cdots \zeta(k)$, and that Theorem \ref{IkeTau} can be
  applied to get the result.
\end{proof}

Now we point out a weaker version of lemma \ref{HNFcount} whose proof is even more elementary.

\begin{lemma}\label{easybound} If $T>0$ then
 $
        \# \{ H \subseteq \Z^k \mid [\Z^k:H] \le T \} \le T^k.
 $
\end{lemma}

\begin{proof}
 As above, we obtain the following formula.
 \begin{equation*}
  \begin{split}
        \# \{ H \subseteq \Z^k \mid [\Z^k:H] \le T \}
        & = \# \{ A \in \HNF \mid \det A \le T \} \\
        & = \sum_{{n_1>0, \dots, n_k>0} \atop {n_1 \cdot \dots \cdot n_k \le T}}
                        n_1^{k-1}  n_2^{k-2}  \cdots  n_{k-1}^1 n_k^0
  \end{split}
 \end{equation*}
 We use it to prove the desired inequality by induction on $k$, the case $k=0$ being clear.
 \begin{equation*}
  \begin{split}
        \# \{ H \subseteq \Z^k \mid [\Z^k:H] \le T \}
        & = \sum_{n_1 = 1}^{\lfloor T \rfloor} n_1^{k-1}  
                  \sum_{{n_2>0, \dots, n_k>0} \atop {n_2 \cdots n_k \le T/n_1}}
                        n_2^{k-2} \cdots  n_{k-1}^1 n_k^0 \\
        & = \sum_{n_1 = 1}^{\lfloor T \rfloor} n_1^{k-1}  
                   \cdot \# \{ H \subseteq \Z^{k-1} \mid [\Z^{k-1}:H] \le T/n_1 \} \\
        & \le \sum_{n_1 = 1}^{\lfloor T \rfloor} n_1^{k-1}  (T / n_1)^{k-1} 
                        \qquad \text{[by induction on $k$]} \\
        & = \sum_{n_1 = 1}^{\lfloor T \rfloor} T^{k-1}
          = {\lfloor T \rfloor} \cdot T^{k-1} \le T^k
  \end{split}
 \end{equation*}
\end{proof}

\section{Volumes}

Recall that a bounded subset $U$ of Euclidean space $\R^k$ is said to have {\sl Jordan content} if
its volume can be approximated arbitrarily well by unions of boxes contained in it or by unions of
boxes containing it, or in other words, that the the characteristic function $\chi_U$ is Riemann
integrable.  Equivalently, the boundary $\partial U$ of $U$ has (Lebesgue) measure zero (see
\cite[Theorem 105.2, Lemma 105.2, and the discussion above it]{MR53:8338}).  If $U$ is a possibly
unbounded subset of $\R^k$ whose boundary has measure zero, its intersection with any ball will
have Jordan content.

Now let's consider the Lie group $G = \Sl_k\R$ as a subspace of the Euclidean
space $M_k\R$ of $k\times k$ matrices.  
Siegel defines a Haar measure on $G$
as follows (see page 341 of \cite{siegel45}).
Let $E$ be a subset of $G$.  Letting $I = [0,1]$ be the unit
interval and considering a number $T > 0$, we may consider the following cones.
\begin{equation*}
 \begin{split}
 I \cdot E         &= \{ t \cdot B \mid B \in E, 0 \le t \le 1\} \\
 T \cdot I \cdot E &= \{ t \cdot B \mid B \in E, 0 \le t \le T\} \\
 \R^+  \cdot E &= \{ t \cdot B \mid B \in E, 0 \le t\}
 \end{split}
\end{equation*}
Observe that if $B \in T \cdot I \cdot E$, then $0 \le \det B \le T^k$.

\begin{definition}\label{muinfdef}
We say that $E$ is measurable if $I \cdot E$ is, and in that case we define
$\mu_\infty(E) = \vol(I \cdot E) \in [0,\infty]$.  
\end{definition}

The Jacobian of left or right
multiplication by a matrix $\gamma$ on $M_k\R$ is $(\det B)^k$, so for
$\gamma \in \Sl_k\R$ volume is preserved.  Thus the measure is invariant
under $G$, by multiplication on either side.  According to Siegel, the
introduction of such invariant measures on Lie groups goes back to Hurwitz
(see \cite[p.~546]{MR27:4723b} or \cite{Hurwitz97}).

Let $F \subseteq G$ be the fundamental domain for the action of $\Gamma = \Sl_k\Z$ on the right of
$G$ presented in \cite[section 7]{MR21:1966}; it's an elementary construction of a fundamental
domain which is a Borel set without resorting to Minkowski's reduction theory.  In each orbit they
choose the element which is closest to the identity matrix in the standard Euclidean norm on $M_k\R
\cong \R^{k^2}$, and ties are broken by ordering $M_k\R$ lexicographically.  This set $F$ is the
union of an open subset of $G$ (consisting of those matrices with no ties) and a countable number
of sets of measure zero.  

The intersection of $T \cdot I \cdot F$ with a ball has Jordan content.  To establish that, it is
enough to show that the measure of the boundary $\partial F$ in $G$ is zero.  Suppose $g \in
\partial F$.  Then it is a limit of points $g_i \not\in F$, each of which has another point $g_i
h_i$ in its orbit which is at least as close to $1$.  Here $h_i$ is in $\Sl_k (\Z)$ and is not $1$.
The sequence $i \mapsto g_i h_i$ is bounded, and thus so is the sequence $h_i$; since $\Sl_k (\Z)$
is discrete, that implies that $h_i$ takes only a finite number of values.  So we may assume $h_i =
h$ is independent of $i$, and is not $1$.  By continuity, $gh$ is at least as close to $1$ as $g$
is.  Now $g$ is also a limit of points $f_i$ in $F$, each of which has $f_i h$ not closer to $1$
than $f_i$ is.  Hence $gh$ is not closer to $1$ than $g$ is, by continuity.  Combining, we see that
$gh$ and $g$ are equidistant from $1$.  The locus of points $g$ in $\Sl_k (\R)$ such that $gh$ and
$g$ are equidistant from $1$ is given by the vanishing of a nonzero quadratic polynomial, hence has
measure zero.  The boundary $\partial F$ is contained in a countable number of such sets, because
$\Sl_k (\Z)$ is countable, hence has measure zero, too.

We remark that $\HNF$ contains a unique representative for each orbit of the
action of $\Sl_k\Z$ on $\{ A \in M_k\Z \mid \det A > 0 \}$.  The same is true
for $\R^+ \cdot F$.  Restricting our attention to matrices $B$ with $\det B \le
T^k$ we see that $ \# ( T \cdot I \cdot F \cap M_k \Z ) = \# \{ A \in \HNF \mid \det A \le T^k \}$.

Warning: $\HNF$ is not contained in $\R^+ \cdot F$.  To convince yourself of
this, consider the matrix $A = \begin{pmatrix}5&-8\\3&5 \end{pmatrix}$ of
determinant $49$.  Column operations with integer coefficients reduce it to
$B = \begin{pmatrix}49&18\\0&1\end{pmatrix}$, but $(1/7)A$ is closer to the
identity matrix than $(1/7)B$ is, so $B \in HNF$, but $B \not\in \R^+ \cdot
F$.

We want to approximate the volume of $T \cdot I \cdot F$ by counting the
lattice points it contains, i.e., by using the number $\# ( T \cdot I \cdot F
\cap M_k \Z )$, at least when $T$ is large.  Alternatively, we may use $\# (
I \cdot F \cap r \cdot M_k \Z )$, when $r$ is small.

\begin{definition}\label{muZdef}
  Suppose $U$ is a subset of $\R^n$.  Let \[N_r(U) = r^n \cdot \# \{ U \cap r
  \cdot \Z^n \}\] and let \[\mu_\Z(U) = \lim_{r \to 0} N_r(U),\] if
  the limit exists, possibly equal to $+\infty$.  An equation involving
  $\mu_\Z(U)$ is to be regarded as true only if the limit exists.
\end{definition}

\begin{lemma}\label{muIF}
  $\mu_\Z(I \cdot F) = \zeta(2)  \zeta(3)  \cdots  \zeta(k) / k$
\end{lemma}

\begin{proof}
  We replace $r$ above with $1/T$:
  \begin{equation*}
   \begin{split}
          \mu_\Z(I \cdot F)
          & = \lim_{T \to \infty} T^{-k^2} \cdot \# ( T \cdot I \cdot F \cap M_k \Z ) \\
          & = \lim_{T \to \infty} T^{-k^2} \cdot \# \{ A \in \HNF \mid \det A \le T^k \} \\
          & = \lim_{T \to \infty} T^{-k^2} \cdot \# \{ H \subseteq \Z^k \mid [\Z^k:H] \le T \} \\
          & = \zeta(2)  \zeta(3)  \cdots  \zeta(k) / k  
                  \qquad \text{[using lemma \ref{HNFcount}]}
   \end{split}
  \end{equation*}
\end{proof}

\begin{lemma}\label{jcon}
  If $U$ is a bounded subset of $\R^n$ with Jordan content, then $\mu_\Z(U) =
  \vol U $.
\end{lemma}

\begin{proof}
  Subdivide $\R^n$ into cubes of width $r$ (and of volume $r^n$) centered at the points of $r\Z^n$.
  The number $\# \{ U \cap r \cdot \Z^n \}$ lies between the number of cubes contained in $U$ and
  the number of cubes meeting $U$, so $r^n \cdot \# \{ U \cap r \cdot \Z^n \}$ is captured between
  the total volume of the cubes contained in $U$ and the total volume of the cubes meeting $U$,
  hence approaches the same limit those two quantities do, namely $\vol U$.
\end{proof}

\begin{lemma}\label{cusps} 
  Let $B_R$ be the ball of radius $R>0$ centered at the origin, and let $U$ be
  a subset of $\R^n$ whose boundary has measure zero.
  \begin{enumerate}
    \item\label{two} For all $R$, the quantity $\mu_\Z(U)$ exists if and only
      if $\mu_\Z(U - B_R)$ exists, and in that case, $\mu_\Z(U) = \vol(U \cap
      B_R) + \mu_\Z(U-B_R)$.
    \item\label{twoo} If $\mu_\Z(U)$ exists then $\mu_\Z(U) = \vol(U) + \lim_{R
        \to \infty} \mu_\Z(U-B_R)$.
    \item\label{one} If $\vol(U) = +\infty$, then $\mu_\Z(U) = +\infty$.
    \item\label{four} If $\lim_{R \to \infty} \limsup_{r \to 0} N_r(U-B_R) =
      0$, then $\mu_\Z(U) = \vol(U)$.
  \end{enumerate}
\end{lemma}

\begin{proof}
  Writing $U = (U \cap B_R) \cup (U - B_R)$ we have 
        $$N_r(U) = N_r(U \cap B_R) + N_r(U - B_R).$$
  For each $R>0$, the set $U \cap B_R$ is a bounded
  set with Jordan content, and thus lemma \ref{jcon} applies to it.  We deduce
  that 
  $$\liminf_{r \to 0} N_r(U) = \vol(U \cap B_R) + \liminf_{r \to 0} N_r(U - B_R)$$
  and
  $$\limsup_{r \to 0} N_r(U) = \vol(U \cap B_R) + \limsup_{r \to 0} N_r(U - B_R),$$
  from which we can deduce (\ref{two}), because $\vol(U \cap B_R) < \infty$.
  We deduce (\ref{twoo}) from (\ref{two}) by taking limits.
  Letting $R \to \infty$ in the equalities above we see that
  $$\liminf_{r \to 0} N_r(U) = \vol(U) + \lim_{R \to \infty} \liminf_{r \to 0} N_r(U - B_R)$$
  and
  $$\limsup_{r \to 0} N_r(U) = \vol(U) + \lim_{R \to \infty} \limsup_{r \to 0} N_r(U - B_R),$$
  in which some of the terms might be $+\infty$.  Now
  (\ref{one}) follows from $\liminf_{r \to 0} N_r(U) \ge \vol(U)$, and
  (\ref{four}) follows because if $$\lim_{R \to \infty} \limsup_{r \to 0}
  N_r(U-B_R) = 0,$$ then $$\lim_{R \to \infty} \liminf_{r \to 0} N_r(U-B_R) =
  0$$ also.
\end{proof}

\begin{lemma} If $U$ is a subset of $\R^n$ whose boundary has measure zero,
  and $ \mu_\Z(U) = \vol(U)$, then $\vol(T \cdot U) \sim \# ( T \cdot U \cap
  \Z^n )$ as $T \to \infty$.
\end{lemma}

\begin{proof}
The statement follows immediately from the definitions.
\end{proof}

Care is required in trying to compute the volume of $I \cdot F$ by counting lattice points in it,
for it is not a bounded set (even for $k=2$, because $\begin{pmatrix}a&0\\0&1/a\end{pmatrix} \in
F$).

\begin{example}\label{unboundedexample}\rm
  It's easy to construct an unbounded region where counting lattice points
  does not determine the volume, by concentrating infinitely many very thin
  spikes along rays of rational slope with small numerator and denominator.
  Consider, for example, a bounded region $B$ in $\R^2$ with Jordan content
  and nonzero area $v = \vol B$, for which (by Lemma \ref{jcon}) $\mu_\Z B =
  \vol B$.  Start by replacing $B$ by its intersection $B'$ with the lines
  through the origin of rational (or infinite) slope -- this doesn't change
  the value of $\mu_\Z$, because every lattice point is contained in a line
  of rational slope, but now the boundary $\partial B'$ does not have measure
  zero.  To repair that, we enumerate the lines $M_1, M_2, \dots$ through the
  origin of rational slope, and for each $i = 1,2,3,\dots$ we replace $R_i =
  B \cap M_i$ by a suitably scaled and rotated version $L_i$ of it contained
  in the line $N_i$ of slope $i$ through the origin, with scaling factor
  chosen precisely so $L_i$ intersects each $r \cdot \Z^2$ in the same number
  of points as does $R_i$, for every $r>0$.  The scaling factor is the ratio
  of the lengths of the shortest lattice points in the lines $M_i$ and $N_i$.
  The union $L = \bigcup L_i$ has $\mu_\Z L = \mu_\Z B = v \ne 0$, but it and
  its boundary have measure zero.
%
%
%
%
\end{example}

\section{Reduction Theory}\label{sec:redthy}

In this section we apply reduction theory to show that the volume of $I \cdot
F$ can be computed by counting lattice points.

We introduce a few basic notions about lattices.  For a more leisurely
introduction see \cite{MR86h:22018}.

\begin{definition}
A {\sl lattice} is a free abelian group $L$ of finite rank
equipped with an inner product on the vector space $L \otimes \R$.
\end{definition}

We will regard $\Z^k$ or one of its subgroups as a lattice by endowing it
with the standard inner product on $\R^k$.

\def\min{\operatorname{min}}
\def\covol{\operatorname{covol}}

\begin{definition}
  If $L$ is a lattice, then a sublattice $L' \subseteq L$ is a subgroup with
  the induced inner product.  The quotient $L/L'$, if it's torsion free, is
  made into a lattice by equipping it with the inner product on the
  orthogonal complement of $L'$.
\end{definition}

There's a way to handle lattices with torsion, but we won't need them.

\begin{definition}
  If $L$ is a lattice, then $\covol L$ denotes the volume of a fundamental
  domain for $L$ acting on $L \otimes \R$.
\end{definition}

The covolume can be computed as $|\det(\theta v_1,\cdots,\theta v_k)|$, where
$\theta : L \otimes \R \to \R^k$ is an isometry,
$\{v_1,\dots,v_k\}$ is a basis of $L$, and $(\theta v_1,\dots,\theta v_k)$ denotes the
matrix whose $i$-th column is $\theta v_i$.  We have the identity $\covol(L) = \covol(L') \cdot
\covol(L/L')$ when $L/L'$ is torsion free.

If $L$ is a subgroup of $\Z^k$ of finite index, then $\covol L = [\Z^k : L]$.

\begin{definition}
If $L$ is a nonzero lattice, then $\min L$ denotes the smallest length of a nonzero
vector in $L$.
\end{definition}

If $L$ is a lattice of rank 1, then $\min L = \covol L$.

\begin{proposition}\label{inequ}
For any natural number $k > 0$, there is a constant $c$ such that for any $S \ge 1$ and for
any $T > 0$ the following inequality holds.
$$
c S^{-k} T^{k^2}
\ge
\# \{ L \subseteq \Z^k \mid [\Z^k: L] \le T^k {\rm \ and\ } \min L \le T/S \}.
$$
\end{proposition}

\begin{proof}
  For $k=1$ we may take $c=2$, so assume $k \ge 2$.  Letting $N$
  be the number of these lattices $L$, we bound $N$ by picking within each $L$ a nonzero vector $v$
  of minimal length, and counting the pairs $(v,L)$ instead.  For each $v$ occurring in such pair
  we write $v$ in the form $v = n_1 v_1$ where $n_1 \in \N$ and $v_1$ is a primitive vector of
  $\Z^k$, and then we extend $\{v_1\}$ to a basis $B = \{v_1,\dots,v_k\}$ of $\Z^k$.  We count the
  lattices $L$ occurring in such pairs with $v$ by putting a basis $C$ for $L$ into Hermite normal
  form with respect to $B$, i.e., it will have the form $C = \{n_1 v_1, A_{12} v_1 + n_2 v_2, \dots
  , A_{1k} v_1 + \dots + A_{k-1,k} v_{k-1} + n_k v_k\}$, with $n_i > 0$ and $0 \le A_{ij} < n_i$.
  Notice that $n_1$ has been determined in the previous step by the choice of $v$.  The number of
  vectors $v \in \Z^k$ satisfying $\| v \| \le T/S$ is bounded by a number of the form $c(T/S)^k$;
  for $c$ we may take a large enough multiple of the volume of the unit ball.  With notation as above, and
  counting the bases for $C$ in Hermite normal form as before, we see that
\begin{equation*}
 \begin{split}
        N & \le \sum_{\| v \| \le T/S} \sum_{{n_2>0, \dots, n_k>0}\atop{n_1 \cdots n_k \le T^k}}
                        n_1^{k-1}  
                        n_2^{k-2}  \cdots 
                        n_{k-1}^{1} 
                        n_k^{0} \\
        & = \sum_{\| v \| \le T/S} n_1^{k-1} \sum_{{n_2>0, \dots, n_k>0}\atop{n_2 \cdots n_k \le T^k/n_1}}
                        n_2^{k-2}  \cdots n_{k-1}^{1} n_k^{0} \\
        & = \sum_{\| v \| \le T/S} n_1^{k-1} \cdot 
                \# \{ H \subseteq \Z^{k-1} \mid [\Z^{k-1}:H] \le T^k/n_1 \} \\
        & \le \sum_{\| v \| \le T/S} n_1^{k-1} (T^k / n_1)^{k-1} 
                        \qquad \text{[by Lemma \ref{easybound}]} \\
        & = \sum_{\| v \| \le T/S}  T^{ k (k-1) } \\
        & \le c(T/S)^k T^{ k (k-1) } \\
        & = c S^{-k} T^{k^2}.
 \end{split}
\end{equation*}
\end{proof}

\begin{corollary}\label{limz}
  The following equality holds.
  \begin{equation*}
    0 = \lim_{S \to \infty}
    \limsup_{T \to \infty}
    T^{-k^2} \cdot \# \{ L \subseteq \Z^k \mid [\Z^k: L] \le T^k {\rm \ and\ } \min L \le T/S \} 
  \end{equation*}
\end{corollary}

The following two lemmas are standard facts.  Compare them, for
example, with \cite[1.4 and 1.5]{MR39:5577}.

\begin{lemma}\label{lifting}
  Let $L$ be a lattice and let $v \in L$ be a primitive vector.  Let $\bar L = L / \Z v$, let
  $\bar w \in \bar L$ be any vector, and let $w \in L$ be a vector of minimal length among all
  those that project to $\bar w$.  Then $\| w \|^2 \le \| \bar w \| ^2 + (1/4) \| v \| ^2$.
\end{lemma}

\begin{proof}
  The vectors $w$ and $w \pm v$ project to $\bar w$, so $\| w \| ^2 \le \|w\pm v\|^2 =
  \|w\|^2+\|v\|^2\pm 2 {\langle}w,v{\rangle}$, and thus $| \langle w,v \rangle | \le (1/2) \|v\|^2$.  We see then that
\[
\begin{split}
  \|\bar w\|^2 & = \| w - \frac{\langle w,v \rangle}{\|v\|^2} v \|^2 \\
               & = \|w \|^2 - \frac{\langle w,v\rangle^2}{\|v\|^2} \\
               & \ge \|w\|^2 - \frac14 \|v\|^2.
\end{split}
\]
\end{proof}

\begin{lemma}\label{plane}
  Let $L$ be a lattice of rank $2$ with a nonzero vector $v \in L$ of minimal length.  Let $L' =
  \Z v$ and $L'' = L/L'$.  Then $\covol L'' \ge (\sqrt 3 / 2) \covol L'$.
\end{lemma}

\begin{proof}
  Let $\bar w \in L''$ be a nonzero vector of minimal length, and lift it to a vector $w \in L$ of
  minimal length among possible liftings.  By lemma \ref{lifting} $\| w \|^2 \le \| \bar w \| ^2 +
  (1/4) \| v \| ^2$.  Combining that with $\|v\|^2 \le \|w\|^2$ we deduce that $\covol L'' = \| \bar w \|
  \ge (\sqrt 3 / 2) \| v \| = (\sqrt 3 / 2) \covol L'$.
\end{proof}

\def\minbasis{\operatorname{minbasis}}

\begin{definition}
If $L$ is a lattice, then $\minbasis L$ denotes the smallest value possible
for $(\| v_1 \|^2 + \dots + \| v_k \|^2)^{1/2}$, where $\{v_1,\dots,v_k\}$ is
a basis of $L$.
\end{definition}

\begin{proposition}\label{minbound}
  Given $k \in \N$ and $S \ge 1$, for all $R \gg 0$, for all $T > 0$, and for all lattices $L$ of
  rank $k$ with $\covol L \le T^k$, if $\minbasis L \ge RT$ then $\min L \le T/S$.
\end{proposition}

\begin{proof}
  We show instead the contrapositive: provided $\covol L \le T^k$, if $\min L > T/S$ then
  $\minbasis L < RT$.  There is an obvious procedure for producing an economical basis of a lattice
  $L$, namely: we let $v_1$ be a nonzero vector in $L$ of minimal length; we let $v_2$ be a vector
  in $L$ of minimal length among those projecting onto a nonzero vector in $L/(\Z v_1)$ of minimal
  length; we let $v_3$ be a vector in $L$ of minimal length among those projecting onto a vector in
  $L/(\Z v_1)$ of minimal length among those projecting onto a nonzero vector in $L/(\Z v_1 + \Z
  v_2)$ of minimal length; and so on.  A vector of minimal length is primitive, so one can show by
  induction that the quotient group $L/(\Z v_1 + \dots + \Z v_i)$ is torsion free; the case where
  $i=k$ tells us that $L = \Z v_1 + \dots + \Z v_k$.  Let $L_i = \Z v_1 + \dots + \Z v_i$, and let
  $\alpha_i = \covol(L_i/L_{i-1})$, so that $\alpha_1 = \| v_1 \| = \min L > T/S$.
  
  Applying Lemma \ref{plane} to the rank $2$ lattice $L_i / L_{i-2}$ shows that $\alpha_{i} \ge A
  \alpha_{i-1}$, where $A = \sqrt 3/2$, and repeated application of Lemma \ref{lifting} shows that
  $\| v_i \|^2 \le \alpha_i^2 + (1/4)( \alpha_{i-1}^2 + \dots + \alpha_{1}^2)$, so of course $
   \| v_i \|^2 \le (1/4)( \alpha_{k}^2 + \dots + \alpha_{i+1}^2) + \alpha_i^2 + (1/4)( \alpha_{i-1}^2
  + \dots + \alpha_{1}^2)$.  We deduce that
  \begin{equation}\label{foo}
     \minbasis L
           \le ( \sum_{i=1}^k \| v_i \|^2 )^{1/2}
           \le \biggl(\frac{k+3}{4} \sum \alpha_i^2\biggr)^{1/2}.
  \end{equation}
  Going a bit further, we see that
  \begin{equation*}
   \begin{split}
          T^k & \ge \covol L \\
              & = \alpha_1 \cdots \alpha_k \\
              & \ge A^{0+1+2+\dots+(i-2)} \alpha_1^{i-1} \cdot A^{0+1+2+\dots+(k-i)} \alpha_{i}^{k-i+1} \\
              & > c_1 (T/S)^{i-1}  \alpha_{i}^{k-i+1}
   \end{split}
  \end{equation*}
  where $c_1$ is some constant depending on $S$ 
  which we may take to be independent of $i$.  Dividing through by
  $T^{i-1}$ we get $T^{k-i+1} > c_2 \alpha_{i}^{k-i+1}$, from which we deduce that $T > c_3
  \alpha_{i}$, where $c_2$ and $c_3$ are new constants (depending only on $S$).  Combining these latter
  inequalities for each $i$, we find that $(((k+3)/4) \sum \alpha_i^2)^{1/2} < R T$, where $R$ is a
  new constant (depending only on $S$); combining that with (\ref{foo}) yields the result.
\end{proof}

\begin{corollary}\label{minbb}
  The following equality holds.
  \[
  0 = \lim_{R \to \infty} \limsup_{T \to \infty} T^{-k^2} \cdot
  \# \{ L \subseteq \Z^k \mid [\Z^k : L] \le T^k {\rm \ and\ } \minbasis L \ge
  RT \}
  \]
\end{corollary}

\begin{proof}
  Combine (\ref{limz}) and (\ref{minbound}).
\end{proof}

If in the definition of our fundamental domain $F$ we had taken the smallest element of each orbit,
rather than the one nearest to $1$, we would have been almost done now.  The next lemma takes care
of that discrepancy.

\def\size{\operatorname{size}}

\begin{definition}
If $L$ is a (discrete) lattice of rank $k$ in $\R^k$, then $\size L$ denotes the value of
$(\| w_1 \|^2 + \dots + \| w_k \|^2)^{1/2}$, where $\{w_1,\dots,w_k\}$ is
the (unique) basis of $L$ satisfying $(w_1,\dots,w_k) \in \R^+ \cdot F$.
\end{definition}

\begin{lemma}\label{minbsize}
For any (discrete) lattice $L \subseteq \R^k$ of rank $k$ the inequalities
 $$\minbasis L \le \size L \le \minbasis L + 2 \sqrt k (\covol L)^{1/k}$$ hold.
\end{lemma}

\begin{proof}
Let $\{v_1,\dots,v_k\}$ be the basis envisaged in the definition of $\minbasis
L$, let $\{w_1,\dots,w_k\}$ be the basis of $L$ envisaged the definition of
$\size L$, and let $U = (\covol L)^{1/k} = (\det (v_1,\dots,v_k))^{1/k} = (\det
(w_1,\dots,w_k))^{1/k}$.  The following chain of inequalities gives the result.
\begin{align}
     \minbasis L & = \|(v_1 , \dots , v_k ) \| \le \size L \notag\\
        & = \|(w_1 , \dots , w_k )\| \le \|(w_1 , \dots , w_k ) - U \cdot 1_k \| + U \sqrt k \notag\\
        & \le \|(v_1 , \dots , v_k ) - U \cdot 1_k \| + U \sqrt k \notag\\
        & \le \|(v_1 , \dots , v_k ) \| + 2 U \sqrt k = \minbasis L + 2 U \sqrt k \notag
\end{align}
\end{proof}

\begin{corollary}\label{sizeb}
  The following equality holds.
  $$
  0 = \lim_{Q \to \infty} \limsup_{T \to \infty} T^{-k^2}
  \cdot \# \{ L \subseteq \Z^k \mid [\Z^k : L] \le T^k {\rm \ and\ } \size L
  \ge Q T \}
  $$
\end{corollary}

\begin{proof}
  It follows from (\ref{minbsize}) that given $R > 0$, for all $Q \gg 0$ (namely $Q \ge R + 2 \sqrt
  k$) if $\covol L \le T^k$ and $\size L \ge QT$ then $\minbasis L \ge RT$.  Now apply
  (\ref{minbb}).
\end{proof}

\begin{theorem}\label{red}
  $ \vol(I \cdot F) = \mu_\Z(I \cdot F) $.
\end{theorem}

\begin{proof}
  Observe that $\# \{ L \subseteq \Z^k \mid \covol L \le T^k {\rm \ and\ } \size L \ge Q T \} = \#
  ((T \cdot I \cdot F - B_{QT}) \cap M_k \Z) = \# ((I \cdot F - B_{Q}) \cap T^{-1} M_k \Z) $, so
  replacing $1/T$ by $r$, Corollary \ref{sizeb} implies that $\lim_{Q \to \infty} \limsup_{r \to
    0} N_r(I \cdot F - B_Q) = 0$, which allows us to apply Lemma \ref{cusps} (\ref{four}).
\end{proof}

The theorem allows us to compute the volume of $F$ arithmetically, simultaneously showing it's finite.

\begin{corollary}
\label{volzeta}
 $ \mu_\infty(G/\Gamma) = \zeta(2)  \zeta(3)  \cdots \zeta(k) / k $
\end{corollary}


\begin{proof}
  Combine the theorem with lemma \ref{muIF} as follows. 
  \begin{equation*}
    \mu_\infty(G/\Gamma) = \mu_\infty(F) = \vol(I \cdot F) =
    \mu_\Z(I \cdot F) = \zeta(2) \zeta(3) \cdots \zeta(k) / k
  \end{equation*}
\end{proof}

\section{$p$-adic volumes}\label{padicvol}

In this section we reformulate the computation of the volume of $G/\Gamma$ to
yield a natural and informative computation of the Tamagawa number of
$\Sl_k$.  We are interested in the form of the proof, not its length, so we
incorporate the proofs of (\ref{volzeta}) and (\ref{zetalatt})
rather than their statements.  The
standard source for information about $p$-adic measures and Tamagawa measures
is Chapter II of \cite{MR83m:10032}, and the proof we simplify occurs there
in sections 3.1 through 3.4.
See also \cite{MR36:171} and \cite{MR31:2249}.

We let $\mu_p$ denote the standard translation invariant measure on $\Q_p$
normalized so that $\mu_p(\Z_p) = 1$.  Let $\mu_p$ also denote the product measure on
the ring of $k$ by $k$ matrices, $M_k(\Q_p)$.  Observe that $\mu_p(M_k(\Z_p)) =
1$.  

For $x \in \Q_p$, let $|x|_p$ denote the standard valuation normalized so
that $|p|_p= 1/p$

If $A \in M_k(\Q_p)$ and $U \subseteq \Q_p^k$, then $\mu_p(A\cdot U) = |\det
A|_p \cdot \mu_p(U)$.  (To prove this, first diagonalize $A$ using row and
column operations, and then assume that $U$ is a cube.)  It follows that if
$V \subseteq M_k(\Q_p)$, then $\mu_p(A\cdot V) = |\det A|_p^k \cdot \mu_p(V)$.

Consider $\Gl_k(\Z_p)$ as an open subset of $M_k(\Z_p)$.
The following computation occurs on page 31 of \cite{MR83m:10032}.
\begin{equation}\label{measGL}
 \begin{split}
   \mu_p(\Gl_k(\Z_p) )
        & = \# (\Gl_k(\F_p)) / p^{k^2} \\
        & = (p^k-1)(p^k-p)\cdots (p^k-p^{k-1}) / p^{k^2} \\
        & = (1 - p^{-k}) (1 - p^{-k+1}) \cdots (1 - p^{-1})
 \end{split}
\end{equation}

Weil considers the open set $M_k(\Z_p)^* = \{A \in M_k(\Z_p) \mid \det A
\ne 0 \}$.

\begin{lemma}\label{Mstar}
 $ \mu_p(M_k(\Z_p)^*) = 1 $
\end{lemma}

\begin{proof}
Let $Z = M_k(\Z_p) \setminus M_k(\Z_p)^*$ be the set of singular matrices.
If $A \in Z$, then one of the columns of $A$ is a linear combination of the
others.  (This depends on $\Z_p$ being a discrete valuation ring -- take any
linear dependency with coefficients in $\Q_p$ and multiply the coefficients
by a suitable power of $p$ to put all of them in $\Z_p$, with at least one of
them being invertible.)  For each $n\ge 0$ we can get an upper bound for the
number of equivalence classes of elements of $Z$ modulo $p^n$ by enumerating
the possibly dependent columns, the possible vectors in the other columns,
and the possible coefficients in the linear combination: $\mu_p(Z) \le
  \lim_{n \to \infty} k \cdot (p^{nk})^{k-1} \cdot (p^n)^{k-1} / (p^n)^{k^2} 
= \lim_{n \to \infty} k \cdot p^{-n}
= 0$.
\end{proof}

We call rank $k$ submodules $J$ of $\Z_p^k$ {\sl lattices}.  To each $A \in
M_k(\Z_p)^*$ we associate the lattice $J = A \Z_p^k \subseteq \Z_p^k$.  This
sets up a bijection between the lattices $J$ and the orbits of $\Gl_k(\Z_p)$
acting on $M_k(\Z_p)^*$.  The measure of the orbit corresponding to $J$ is
$\mu_p(A \cdot \Gl_k(\Z_p)) = |{\det A}|_p^k \cdot \mu_p(\Gl_k(\Z_p)) =
[\Z_p^k:J]^{-k} \cdot \mu_p(\Gl_k(\Z_p))$.  Now we sum over the orbits.
\begin{equation}\label{localTam}
 \begin{split}
        1
        & = \mu_p (M_k(\Z_p)^*) \\
        & = \sum_J \Bigl([\Z_p^k:J]^{-k} \cdot \mu_p(\Gl_k(\Z_p))\Bigr) \\
        & = \Bigl(\sum_J [\Z_p^k:J]^{-k}\Bigr) \cdot \mu_p(\Gl_k(\Z_p))
 \end{split}
\end{equation}
An alternative way to prove (\ref{localTam}) would be to use the local
analogue of (\ref{zetalatt}), which holds and asserts that $\sum_J
[\Z_p^k:J]^{-s} = (1-p^{k-1-s})^{-1} (1-p^{k-2-s})^{-1} \cdots
(1-p^{-s})^{-1}$; we could substitute $k$ for $s$ and compare with the number
in (\ref{measGL}).  The approach via lemma \ref{Mstar} and (\ref{localTam})
is preferable because $M_k(\Z_p)^*$ provides natural glue that makes the
computation seem more natural.

The product $\prod_p \mu_p(\Gl_k(\Z_p))$ doesn't converge because $\prod_p
(1-p^{-1})$ doesn't converge, so consider the following formula instead.
\begin{equation*}
        1 = \Bigl((1-p^{-1}) \sum_J [\Z_p^k:J]^{-k}\Bigr) \cdot 
        \Bigl( (1-p^{-1})^{-1} \mu_p(\Gl_k(\Z_p)) \Bigr)
\end{equation*}
Now we can multiply these formulas together.
\begin{equation}\label{globalTam}
        1 = \Bigl(\prod_p (1-p^{-1}) \sum_J [\Z_p^k:J]^{-k}\Bigr) \cdot 
        \prod_p \Bigl((1-p^{-1})^{-1} \mu_p(\Gl_k(\Z_p)) \Bigr)
\end{equation}
We've parenthesized the formula above so it has one factor for each place of
$\Q$, and now we connect each of them with a volume involving $\Sl_k$ at that
place.

We use the Haar measure on $\Sl_k(\Z_p)$ normalized to have total volume $$\#
\Sl_k(\F_p) / p^{\dim \Sl_k} .$$ The normalization anticipates (\ref{meas}),
which shows how a gauge form could be used to construct the measure, or
alternatively, it ensures that the exact sequence $1 \to \Sl_k(\Z_p) \to
\Gl_k(\Z_p) \to  \Z_p^\times \to 1$ of groups leads to the desired
assertion $ \mu_p(\Gl_k(\Z_p)) = \mu_p(\Z_p^\times) \cdot \mu_p(\Sl_k(\Z_p))$
about multiplicativity of measures.  We rewrite the factor of the right hand
side of (\ref{globalTam}) corresponding to the prime $p$ as follows.
\begin{equation}\label{fubini}
\begin{split}
   (1-p^{-1})^{-1} \mu_p(\Gl_k(\Z_p)) 
& = \mu_p(\Z_p^\times)^{-1} \cdot \mu_p(\Gl_k(\Z_p)) \\
& = \mu_p(\Sl_k(\Z_p)).
\end{split}
\end{equation}

To evaluate the left hand factor of the right hand side of (\ref{globalTam}),
we insert the complex variable $s$.  Because the ring $\Z$ is a principal
ideal domain, any finitely generated sub-$\Z$-module $H \subseteq \Z^k$ is
free.
Hence a lattice $H
\subseteq \Z^k$ is determined freely by its localizations $H_p = H \otimes_\Z
\Z_p \subseteq \Z_p^k$ (where $H_p = \Z_p^k$ for all but finitely many $p$),
and its index is given by the formula
\begin{equation}\label{localglobal}
        [\Z^k:H] = \prod_p [\Z_p^k:H_p],
\end{equation}
in which only a finite number of terms are not equal
to~$1$.
\begin{equation}\label{work1}
 \begin{split}
      \res_{s = k} & \zeta(\Z^k,s) \\
  & = \res_{s = k} \sum_{H} [\Z^k:H]^{-s} \qquad \text{[by (\ref{zetadef})]} \\
  & = \lim_{s \to k+} \zeta(s-k+1)^{-1} \cdot  \sum_{H} [\Z^k:H]^{-s} \\
  & = \lim_{s \to k+} \biggl( 
          \zeta(s-k+1)^{-1} \bigl( \sum_{H \subseteq \Z^k} \prod_p [\Z_p^k:H_p]^{-s} \bigr)
          \biggr)
         \qquad \text{[by (\ref{localglobal})]} \\
  & = \lim_{s \to k+} \biggl( 
          \zeta(s-k+1)^{-1} \bigl( \prod_p \sum_{J \subseteq \Z_p^k} [\Z_p^k:J]^{-s} \bigr)
          \biggr)
          \qquad \text{[positive terms]} 
         \\
  & = \lim_{s \to k+} \prod_p \bigl( (1-p^{-s+k-1}) \sum_J [\Z_p^k:J]^{-s} \bigr) 
        \\
  & = \prod_p (1-p^{-1}) \sum_J [\Z_p^k:J]^{-k}
       \\
 \end{split}
\end{equation}
Starting again we get the following chain of equalities.
\begin{equation}\label{work2}
 \begin{split}
      \res_{s = k} \zeta(\Z^k,s)
  & = \res_{s = k} \zeta(s-k+1)  \zeta(s-k+2) \cdots  \zeta(s-1)  \zeta(s) \\
  & = \zeta(2) \cdots  \zeta(k-1)  \zeta(k) \\
  & = k \cdot \lim_{T \to \infty} T^{-k} \# \{ H \subseteq \Z^k \mid [\Z^k:H] \le T \} 
        \qquad \text{[by \ref{HNFcount}]} \\
  & = k \cdot \lim_{T \to \infty} T^{-k^2} \# \{ H \subseteq \Z^k \mid [\Z^k:H] \le T^k \} \\
  & = k \cdot \lim_{T \to \infty} T^{-k^2} \# \{ A \in \HNF \mid \det A \le T^k \} \\
  & = k \cdot \mu_\Z (I \cdot F ) \qquad \text{[by definition \ref{muZdef}]} \\
  & = k \cdot \mu_\infty ( \Sl_k(\R) / \Sl_k(\Z) ) \qquad \text{[by \ref{red} and \ref{muinfdef}]}
 \end{split}
\end{equation}

Combining (\ref{work1}) and (\ref{work2}) we get the following equation.
\begin{equation}\label{work3}
   \prod_p (1-p^{-1}) \sum_J [\Z_p^k:J]^{-k}
   =
   k \cdot \mu_\infty ( \Sl_k(\R) / \Sl_k(\Z) )
\end{equation}

We combine (\ref{globalTam}), (\ref{fubini}) and (\ref{work3})
to obtain the following equation.
\begin{equation}\label{almost}
1 = k \cdot \mu_\infty ( \Sl_k(\R) / \Sl_k(\Z) ) \cdot \prod_p \mu_p(\Sl_k(\Z_p))
\end{equation}

To relate this to the Tamagawa number we have to introduce a gauge form
$\omega$ on the algebraic group $\Sl_k$ over $\Q$, invariant by left
translations, as in sections 2.2.2 and 2.4 of \cite{MR83m:10032}.  We can
even get gauge forms over $\Z$.  Let $X$ be a generic element of $\Gl_k$.
The entries of the matrix $X^{-1} dX$ provide a basis for the 1-forms
invariant by left translation on $\Gl_k$.  On $\Sl_k$ we see that $\tr(X^{-1}
dX) = d(\det X) = 0$, so omitting the element in the $(n,n)$ spot will
provide a basis of the invariant forms on $\Sl_k$.  We let $\omega$ be the
exterior product of these forms.  Just as in the proof of Theorem 2.2.5 in
\cite{MR83m:10032} we obtain the following equality.
\begin{equation}\label{meas}
\int_{\Sl_k(\Z_p)} \omega_{p} = \mu_p(\Sl_k(\Z_p))
\end{equation}
The measure $\omega_p$ is defined in \cite[2.2.1]{MR83m:10032} in a neighborhood of a point $P$ by
writing $\omega = f \, dx_1 \wedge \dots \wedge dx_n$ and setting $\omega_p = |f(P)|_p (dx_1)_p
\dots (dx_n)_p$, where $(dx_i)_p$ is the Haar measure on $\Q_p$ normalized so that $\int_{\Z_p}
(dx_i)_p = 1$, and $|c|_p$ is the $p$-adic valuation normalized so that $d(cx)_p = |c|_p (dx)_p$.

  Now we want to determine the constant that relates our
original Haar measure $\mu_\infty$ on $\Sl_k(\R)$ to the one determined by
$\omega_\infty$.  For this purpose, it will suffice to evaluate both measures
on the infinitesimal parallelepiped $B$ in $\Sl_k(\R)$ centered at the
identity matrix and spanned by the tangent vectors $\varepsilon e_{ij}$ for $i
\ne j$ and $\varepsilon ( e_{ii} - e_{kk} )$ for $i < k$.  Here $\varepsilon$ is an
infinitesimal number, and $e_{ij}$ is the matrix with a $1$ in position
$(i,j)$ and zeroes elsewhere.  For the purpose of this computation, we may even
take $\varepsilon = 1$.  We compute easily that $\int_B \omega_\infty = 1$
and
\begin{equation}\label{normaliz}
\begin{split}
\mu_\infty(B) &= \vol(I \cdot B) \\
&= (1/k^2) \cdot |\det(e_{11}-e_{kk},\cdots,e_{k-1,k-1}-e_{kk},\sum e_{ii})| \\
&= (1/k^2) \cdot |\det(e_{11}-e_{kk},\cdots,e_{k-1,k-1}-e_{kk},k e_{kk})| \\
&= (1/k^2) \cdot |\det(e_{11},\cdots,e_{k-1,k-1},k e_{kk})| \\
&= 1/k
\end{split}
\end{equation}
We obtain the following equation.
\begin{equation}\label{relate}
\mu_\infty ( \Sl_k(\R) / \Sl_k(\Z) ) = \frac 1 k \int_{ \Sl_k(\R) / \Sl_k(\Z) } \omega_\infty
\end{equation}
See \cite[\S 14.12, (3)]{MR99g:20090} for an essentially equivalent proof of
this equation.
We may now rewrite (\ref{almost}) as follows.
\begin{equation}\label{final}
1 = \int_{\Sl_k(\R) / \Sl_k(\Z)} \omega_\infty \cdot \prod_p \int_{\Sl_k(\Z_p)} \omega_{p}
\end{equation}
(If done earlier, this computation would have justified normalizing
$\mu_\infty$ differently.)

The Tamagawa number $\tau(\Sl_{k,\Q}) = \int_{\Sl_k(\A_\Q)/\Sl_k(\Q)} \omega$
is the same as the right hand side of (\ref{final}) because $F \times \prod_p
\Sl_k(\Z_p)$ is a fundamental domain for the action of $\Sl_k(\Q)$ on
$\Sl_k(\A_\Q)$.  Thus $\tau(\Sl_{k,\Q}) = 1$.  This was originally proved by
Weil in Theorem 3.3.1 of \cite{MR83m:10032}.  See also \cite{MR35:4226},
\cite{MR90e:11075}, and \cite[\S 14.11, Corollary to Langlands'
Theorem]{MR99g:20090}.

See also \cite[\S 8]{MR35:2900} for an explanation that Siegel's
measure formula amounts to the first determination that $\tau(SO)=2$.

\appendix

\section{Dirichlet series}

\begin{theorem}\label{convergence}
  Suppose we are given a Dirichlet series $f(s) := \sum_{n=1}^\infty a_n n^{-s}$ with nonnegative
  coefficients.  Let $A(T) := \sum_{n \le T} a_n$.  If $ A(T) = O(T^k)$ as $T \to \infty$, then
  $\sum_{n=T}^\infty a_n n^{-s} = O(T^{k-s})$ as $T \to \infty$, and thus $f(s)$ converges for all
  complex numbers $s$ with  $\re s > k$.
\end{theorem}

\begin{proof}
  Write $\sigma = \re s$ and assume $\sigma > k$.  We estimate the tail of the series as follows.
  \begin{equation*}
    \begin{split}
            \sum_{n=T}^\infty a_n n^{-s}
            & = \int_{T}^\infty x^{-s} \, d A(x) \\
            & = x^{-s} A(x) \big ]_{T}^\infty -   \int_{T}^\infty A(x)  \, d(x^{-s}) \\
            & = x^{-s} A(x) \big ]_{T}^\infty + s \int_{T}^\infty x^{-s-1} A(x)  \, dx \\
            & = O(x^{k-\sigma})  \big ]_T  ^\infty + s \int_{T}^\infty x^{-s-1} O(x^k) \, dx \\
            & = O(T^{k-\sigma}) +                    s \int_{T}^\infty     O(x^{k-\sigma-1}) \, dx \\
            & = O(T^{k-\sigma})
    \end{split}
  \end{equation*}
\end{proof}

\begin{theorem}\label{product}
  Suppose we are given two Dirichlet series 
  \begin{equation*}
      f(s)  := \sum_{n=1}^\infty a_n n^{-s} \qquad \qquad
      g(s)  := \sum_{n=1}^\infty b_n n^{-s}
  \end{equation*}
  with nonnegative coefficients and corresponding coefficient summatory functions
  \begin{equation*}
    A(T) := \sum_{n \le T} a_n 
    \qquad \qquad
    B(T) := \sum_{n \le T} b_n
  \end{equation*}
  Assume that $A(T) = O(T^i)$ and $B(T) = cT^k + O(T^j)$, where $i \le j < k$.
  Let $h(s) := f(s) g(s) = \sum_{n=1}^\infty c_n n^{-s}$,
  and let $C(T) := \sum_{n \le T} c_n$.
  Then $C(T) = c f(k) T^k + O(T^j \log T)$ if $i=j$,
  and $C(T) = c f(k) T^k + O(T^j)$ if $i<j$.
\end{theorem}

\begin{proof}
  The basic idea for this proof was told to us by Harold Diamond.

  Observe that Theorem \ref{convergence} ensures that $f(k)$ converges.  Let's fix the notation
  $\beta(T) = O(\gamma(T))$ to mean that there is a constant $C$ so that $|\beta(T)| \le C
  \gamma(T)$ for all $T \in [1,\infty)$, and simultaneously replace $O(T^j \log T)$ in the
  statement by $O(T^j (1 + \log T))$ in order to avoid the zero of $\log T$ at $T = 1$.  We will
  use the notation in an infinite sum only with a uniform value of the implicit constant $C$.

  We examine $C(T)$ as follows.
  \begin{equation*}
    \begin{split}
      C(T) & = \sum_{n \le T} c_n
             = \sum_{n \le T} \sum_{pq=n} a_p b_q
             = \sum_{pq \le n} a_p b_q \\
           & = \sum_{p \le T} a_p \sum_{q \le T/p} b_q
             = \sum_{p \le T} a_p B(T/p) \\
           & = \sum_{p \le T} a_p \{ c (T/p)^k + O((T/p)^j) \} \\
           & = c T^k \sum_{p \le T} a_p p^{-k} + O(T^j) \sum_{p \le T} a_p p^{-j} \\
           & = c T^k \{ f(k) + O(T^{i-k}) \} + O(T^j) \sum_{p \le T} a_p p^{-j} \\
           & = c f(k) T^k  + O(T^{i}) + O(T^j) \sum_{p \le T} a_p p^{-j}
    \end{split}
  \end{equation*}

  If $i < j$ then $\sum_{p \le T} a_p p^{-j} \le f(j) = O(1)$.
  Alternatively, if $i = j$, then 
  \begin{equation*}
    \begin{split}
       \sum_{p \le T} a_p p^{-j}
       & = \sum_{p \le T} a_p p^{-i}
         = \int_{1-}^T p^{-i} \, d(A(p)) \\
       & = p^{-i} A(p) \bigr]_{1-}^T - \int_{1-}^T A(p) \, d(p^{-i}) \\
       & = T^{-i} A(T) + i \int_{1-}^T A(p) p^{-i-1} \, dp \\
       & = O(1) + O(\int_{1-}^T p^{-1} \, dp)
         = O(1 + \log T)
    \end{split}
  \end{equation*}
  In both cases the result follows.
\end{proof}

The proof of the following ``Abelian'' theorem for generalized Dirichlet series is elementary.

\begin{theorem}\label{elementary}
  Suppose we are given numbers $R$, $k \ge 1$, and $1 \le \lambda_1 \le \lambda_2 \le \dots \to \infty $.
  Suppose that
  \begin{equation*}
    N(T) := \sum_{\lambda_n \le T} 1 = (R + o(1)) \frac{T^k}{k} \qquad ( T \to \infty )
  \end{equation*}
  for some number $R$.  Then the generalized Dirichlet series $\psi(s) := \sum \lambda_n^{-s}$
  converges for all real numbers $s > k$, and $\lim_{s \to k+}(s-k)\psi(s) = R$.
\end{theorem}

\begin{proof}
  In the case $R \ne 0$, the proof can be obtained by adapting the argument
  in the last part of the proof of \cite[Chapter 5, Section 1, Theorem
  3]{MR33:4001}: roughly, one reduces to the case where $k=1$ by a simple
  change of variables, shows $\lambda_n \sim n/R $, uses that to compare a
  tail of $\sum \lambda_n^{-s}$ to a tail of $\zeta(s) = \sum n^{-s}$, and
  then uses $\lim_{s \to 1+}(s-1)\zeta(s) = 1$.
  
  Alternatively, one can refer to \cite[Theorem 10, p. 114]{MR97e:11005b} for
  the statement about convergence, and then to \cite[Theorem 2, p.
  219]{MR97e:11005b} for the statement about the limit.  Actually, those two
  theorems are concerned with Dirichlet series of the form $F(s) = \sum a_n
  n^{-s}$, but the first step there is to consider the growth rate of
  $\sum_{n \le x} a_n$ as $x \to \infty$.  Essentially the same proof works
  for $F(s) = \psi(s)$ by considering the growth rate of $N(x)$ instead.
  
  The result also follows from the following estimate, provided to us by Harold Diamond.  Assume $s
  > k$.
  \begin{equation*}
    \begin{split}
            \psi(s) 
            & := \sum \lambda_n^{-s} \\
            & = \int_{1-}^\infty x^{-s} \, dN(x) \\
            & = x^{-s} N(x) \big ]_{1-}^\infty + s \int_{1}^\infty x^{-s-1} N(x) \, dx \\
            & = O(x^{k-s})  \big ]     ^\infty + s \int_{1}^\infty x^{-s-1} (R + o(1))
                \frac{x^k}{k} \, dx                 \qquad ( x \to \infty ) \\
            & = \frac{s ( R + o(1) )}{k}  
                \int_{1}^\infty x^{-s-1+k} \, dx    \qquad ( s \to k+ ) \\
            & = \frac{s ( R + o(1) )}{k(s-k)}       \qquad ( s \to k+ ) 
    \end{split}
  \end{equation*}
  Notice the shift in the meaning of $o(1)$ from one line to the next, verified by writing
  $\int_{1}^\infty = \int_{1}^b + \int_{b}^\infty$ and letting $b$ go to $\infty$; it turns out that
  for sufficiently small $\epsilon$ the major contribution to $\int_1^\infty x^{-1-\epsilon} \, dx$
  comes from $\int_b^\infty x^{-1-\epsilon} \, dx$.
\end{proof}

The following Wiener-Ikehara ``Tauberian'' theorem is a converse to the previous theorem, but the
proof is much harder.

\begin{theorem}\label{IkeTau}
  Suppose we are given numbers $R>0$,
  $k>0$, $1 \le \lambda_1 \le \lambda_2 \le \dots \to \infty $, and
  nonnegative numbers $a_1, a_2, \dots$.  
  Suppose that the Dirichlet series
  $\psi(s) = \sum a_n \lambda_n^{-s}$ converges for all complex numbers with
  $\Re s > k$, and that the function $\psi(s) - R/(s-k)$ can be extended to a
  function defined and continuous for $\Re s \ge k$.  Then 
  \begin{equation*}
    \sum_{\lambda_n \le T} a_n \sim R T^k / k.
  \end{equation*}
\end{theorem}

\begin{proof}
  Replacing $s$ by $k s$ allows us to reduce to the case where $k=1$, which can
  be deduced directly 
  from the Landau-Ikehara Theorem in \cite{Bochner},
  from Theorem 2.2 on p.~93 of \cite{MR88j:40011},
  from Theorem 1 on p.~464 of \cite{MR50:268},
  or from Theorem 1 on p.~534 of \cite{MR91h:11107}.
  See also Theorem 17 on p.~130 of
  \cite{Wiener33} for the case where $\lambda_n = n$, which suffices for our
  purposes.  
  A weaker prototype of this theorem was first proved by Landau in 1909 \cite[\S 241]{MR16:904d}.
  Other relevant papers include \cite{Wiener32}, \cite{MR16:921e}, and \cite{MR12:405a}.
  See also Bateman's discussion in \cite[Appendix, page 931]{MR16:904d} and
  the good exposition of Abelian and Tauberian theorems in chapter 5 of \cite{MR3:232d}.
\end{proof}

\bibliographystyle{plain} 

\bibliography{papers}
\end{document}